\input ./style/arxiv-vmsta.cfg
\documentclass[numbers,compress,v1.0.1]{vmsta}

\volume{3}
\issue{1}
\pubyear{2016}
\firstpage{59}
\lastpage{78}
\doi{10.15559/16-VMSTA51}

\usepackage{enumerate}

\startlocaldefs

\allowdisplaybreaks

\newtheorem{thm}{Theorem}

\newtheorem{lemma}{Lemma}
\newtheorem{cor}{Corollary}

\theoremstyle{definition}
\newtheorem{defin}{Definition}
\newtheorem{exm}{Example}

\newcommand{\me}{\mathbf}
\newcommand{\mr}{\mathbb}
\newcommand{\md}{\mathcal}
\newcommand{\mt}{\mathsf}

\newcommand{\ld}{\left}
\newcommand{\rd}{\right}

\newcommand{\ip}{\int_{-\pi}^{\pi}}

\endlocaldefs

\begin{document}
\begin{frontmatter}

\title{Minimax interpolation of sequences with stationary increments
and cointegrated sequences}\vspace*{-9pt}

%




\author{\inits{M.}\fnm{Maksym}\snm{Luz}}\email{maksim\_luz@ukr.net}
\author{\inits{M.}\fnm{Mikhail}\snm{Moklyachuk}\corref
{cor1}}\email
{mmp@univ.kiev.ua}
\cortext[cor1]{Corresponding author.}

\address{Taras Shevchenko National University of Kyiv, Kyiv, Ukraine}

\markboth{M. Luz, M. Moklyachuk}{Interpolation of sequences with
stationary increments and cointegrated sequences}

\begin{abstract}
We consider the problem of optimal estimation of the linear
functional $A_N{\xi}=\sum_{k=0}^{N}a(k)\xi(k)$ depending
on the unknown values of a stochastic sequence $\xi(m)$ with
stationary increments from observations of the sequence $\xi(m)+\eta
(m)$ at
points of the set $\mr Z\setminus\{0,1,2,\ldots,N\}$, where $\eta(m)$
is a stationary sequence uncorrelated with $\xi(m)$.
We propose formulas for calculating the mean square error and the spectral
characteristic of the optimal linear estimate of the functional in the
case of spectral certainty, where spectral densities of the sequences
are exactly known. We also consider the problem for a class of
cointegrated sequences.
We propose relations that determine the least favorable spectral
densities and
the minimax spectral characteristics in the case of spectral
uncertainty, where spectral densities are not exactly known while a set
of admissible spectral densities is specified.
\end{abstract}


\received{11 March 2016}
%
\revised{16 March 2016}
\accepted{17 March 2016}
\publishedonline{1 April 2016}
\end{frontmatter}

\section{Introduction}
In this paper, we investigate the problem of estimating the missed
observations of stochastic sequences with stationary increments.
Kolmogorov \cite{Kolmogorov}, Wiener \cite{Wiener}, and Yaglom \cite
{Yaglom:1987a,Yaglom:1987b} developed effective methods of estimation
of the unknown values of stationary sequences and processes.
Later on Yaglom \cite{Yaglom:1955a} and Pinsker
\cite{Pinsk:1955}
introduced and investigated stochastic processes with stationary
increments of order~$n$. Properties of these and
other processes
generalizing the concept of stationarity are described in the books by
Yaglom \cite{Yaglom:1987a,Yaglom:1987b}. The stationary and related
stochastic sequences are widely used in econometrics and
in financial time series analysis. Examples of these sequences are
autoregressive sequences (AR), moving-average sequences (MA), and
autoregressive moving-average sequences (ARMA).
Time series with trends are described by integrated ARMA sequences
(ARIMA) and seasonal time series, which are examples of stochastic
sequences with stationary increments. These models are properly
described in the book by Box, Jenkins, and Reinsel
\cite{Box:Jenkins}. Granger \cite{Granger} introduced a concept of
cointegrated sequences, namely, the integrated sequences such that some
linear combination of them has a lower order of integration.
Cointegrated sequences are described in more details in the paper by
Engle and Granger \cite{Engle:Granger}. We also refer to the papers
\cite{Chigira:Yamamoto,Clements,Gregoir,Johansen} for recent developments.

Traditional methods of finding solutions to extrapolation,
interpolation, and filtering problems for stationary and related
stochastic processes are developed under the basic assumption that the
spectral densities of the considered stochastic processes are exactly known.
However, in most practical situations, complete information on the
spectral densities of the processes is not available.
Investigators can apply the traditional methods considering the
estimated spectral densities instead of the true ones. However, as it
was shown by Vastola and Poor \cite{Vastola} with the help of some
examples, this approach can result in significant increasing of the
value of the error of estimate. Therefore, it is reasonable to derive
estimates that are optimal for all densities from a certain class of
spectral densities. These estimates are called minimax-robust since
they minimize the maximum of the mean-square errors for all spectral
densities from a set of admissible spectral densities simultaneously.
This approach to study the problem of extrapolation of stationary
stochastic processes was introduced by Grenander \cite{Grenander}.
Franke \cite{Franke:1985} investigated the minimax extrapolation and
interpolation problems for stationary sequences applying the convex
optimization methods.
In the book by Moklyachuk \cite{Moklyachuk:2008}, the minimax-robust
estimates of the linear functionals of stationary sequences and
processes are presented. See also the survey paper \cite{Moklyachuk:2015},
The classical and minimax-robust problems of interpolation,
extrapolation, and filtering of the functional of stochastic sequences
with stationary increments are investigated in the papers by Luz and
Moklyachuk \cite{Luz5,Luz8,Luz2,Luz3,Luz6}. Particularly, the
cointegrated sequences are investigated in the papers \cite{Luz5,Luz8}.
The classical extrapolation problem in the case where both the
signal and the noise processes are not stationary was investigated by
Bell \cite{Bell}.

In the present paper, we consider the problem of estimation of the
linear functional\vspace{-3pt}
\[A_N\xi=\sum_{k=0}^{N}a(k)\xi(k),\] which depends on the unknown values
of the sequence $\xi(k)$ with stationary $n$th increments based on
observations of the sequence $\xi(k)+\eta(k)$ at points $m\in\mr
Z\setminus\{0,1,2,\ldots,N\}$. The sequence $\eta(k)$ is assumed to be
stationary and uncorrelated with $\xi(k)$.\looseness=-1

\section{Stationary increment stochastic sequences. Spectral representation}
In this section, we present the main results of the spectral theory of
stochastic sequences with $n$th stationary increments. For more
details, we refer to the books by Yaglom \cite{Yaglom:1987a,Yaglom:1987b}.
\begin{defin}
For a given stochastic sequence $\{\xi(m),m\in\mathbb Z\}$, the
sequence
\begin{equation}
\label{oznachPryrostu}
\xi^{(n)}(m,\mu)=(1-B_{\mu})^n\xi(m)=\sum_{l=0}^n(-1)^l {{n
\choose l}}
\xi(m-l\mu),
\end{equation}
where $B_{\mu}$ is the backward shift operator with step $\mu\in
\mathbb Z$ such that $B_{\mu}\xi(m)=\xi(m-\mu)$, is called
a stochastic $n$th increment sequence with step $\mu\in\mathbb Z$.
\end{defin}

\begin{defin}
The stochastic $n$th increment sequence $\xi^{(n)}(m,\mu)$ generated
by a sto\-chastic sequence $\{\xi(m),m\in\mathbb Z\}$ is wide sense
stationary if the mathematical expectations
\[
\mathsf{E} \xi^{(n)}(m_0,\mu)=c^{(n)}(\mu),
\]
\[
\mathsf{E}
\xi^{(n)}(m_0+m,\mu_1)\overline{\xi^{(n)}(m_0,\mu
_2)}=D^{(n)}(m,\mu
_1,\mu_2)
\]
exist for all $m_0,\mu,m,\mu_1,\mu_2$ and do not depend on $m_0$.
The function $c^{(n)}(\mu)$ is called the mean value of the $n$th
increment sequence, and the function $D^{(n)}(m,\mu_1,\mu_2)$ is
called the structural function of the stationary $n$th increment
sequence (or the structural function of $n$th order of the stochastic
sequence $\{\xi(m),m\in\mathbb Z\}$).
\end{defin}

\begin{thm} \label{thm1}
The mean value $c^{(n)}(\mu)$ and the structural function
$D^{(n)}(m,\mu_1,\mu_2)$ of the stochastic stationary $n$th
increment sequence $\xi^{(n)}(m,\mu)$ can be represented in the following
forms:
\begin{equation}
\label{serFnaR}
c^{(n)}(\mu)=c\mu^n,
\end{equation}
\begin{equation}
\label{strFnaR}
D^{(n)}(m,\mu_1,\mu_2)=\int_{-\pi}^{\pi}e^{i\lambda
m} \bigl(1-e^{-i\mu_1\lambda}\bigr)^n\bigl(1-e^{i\mu_2\lambda
}\bigr)^n\frac{1}
{\lambda^{2n}}dF(\lambda),
\end{equation}
where $c$ is a constant, $F(\lambda)$ is a left-continuous
nondecreasing bounded function with $F(-\pi)=0$. The constant $c$
and the function $F(\lambda)$ are determined uniquely by the
increment sequence $\xi^{(n)}(m,\mu)$.
\end{thm}

Representation (\ref{strFnaR}) and the Karhunen theorem
\cite{Gihman:Skorohod} give us a spectral representation of
the stationary $n$th increment sequence $\xi^{(n)}(m,\mu)$:
\begin{equation}
\label{predZnaR} \xi^{(n)}(m,\mu)=\ip
e^{im\lambda}\bigl(1-e^{-i\mu\lambda}\bigr)^n\frac{1}{(i\lambda
)^n}dZ_{\xi
^{(n)}}(\lambda),
\end{equation}
where $Z_{\xi^{(n)}}(\lambda)$ is a random process with uncorrelated
increments on $[-\pi,\pi)$ with respect to the spectral function
$F(\lambda)$:
\[
\mt E\big|Z_{\xi^{(n)}}(t_2)-Z_{\xi^{(n)}}(t_1)\big
|^2=F(t_2)-F(t_1)\quad
\forall-\pi\leq t_1<t_2<\pi.
\]

We will use the spectral representation (\ref{predZnaR}) for deriving
the optimal linear estimates of unknown values
of stochastic sequences with stationary increments.

\section{Hilbert space projection method of interpolation}
Consider a stochastic sequence $\{\xi(m),m\in\mr Z\}$ with stationary
$n$th increments\break$\xi^{(n)}(m,\mu)$
and uncorrelated with $\xi(m)$ stationary stochastic sequence $\{\eta
(m),m\in\mr Z\}$. Suppose that these sequences have absolutely
continuous spectral functions $F(\lambda)$ and $G(\lambda)$ with
spectral densities
$f(\lambda)$ and $g(\lambda)$, respectively. We will suppose that the
stationary increment $\xi^{(n)}(m,\mu)$ and the stationary sequence
$\eta(m)$ have zero mean values and $\mu>0$.

Interpolation problem for the sequences $\xi(m)$ and $\eta(m)$ is
considered as the problem of the mean-square optimal estimation of the
linear functional
\[
A_N\xi=\sum_{k=0}^{N}a(k)\xi(k),
\]
which depends on the unknown values of the stochastic sequence $\xi(m)$
at points $m=0,1,\ldots,N$ based on observations of the sequence
$\zeta
(m)=\xi(m)+\eta(m)$ at points of the set
$\mr Z\setminus\{0,1,2\ldots,N\}$.

Suppose that the spectral densities $f(\lambda)$ and $g(\lambda)$
satisfy the minimality condition
\begin{equation}
\ip\frac{\lambda^{2n}}{|1-e^{i\lambda\mu}|^{2n}(f(\lambda
)+\lambda
^{2n}g(\lambda))}
d\lambda<\infty.
\label{umova11_i_st.n_d}
\end{equation}
Under this condition, the mean-square error of the estimate of the
functional $A_N\xi$ is not equal to zero \cite{Rozanov}.

The functional $A_N\xi$ admits the representation
\begin{equation}\label{zobrazh A_N_i_st.n_d}
A_N\xi=A_N\zeta-A_N\eta=B_N\zeta-A_N\eta-V_N\zeta=H_N\xi
-V_N\zeta,
\end{equation}
where
\[
H_N\xi:=B_N\zeta-A_N\eta,\qquad A_N\zeta=\sum_{k=0}^{N}a(k)\zeta
(k),\qquad A_N\eta=\sum_{k=0}^{N}a(k)\eta(k),
\]
\[
B_N\zeta=\sum_{k=0}^{N}b_{\mu,N}(k)\zeta^{(n)}(k,\mu),\qquad
V_N\zeta=\sum_{k=-\mu n}^{-1}v_{\mu,N}(k)\zeta(k).
\]
The coefficients $v_{\mu,N}(k)$, $k=-\mu n,-\mu n+1,\ldots,-1$, and
$b_{\mu,N}(k)$, $k=0,1,2,\ldots,\break N$, are calculated by the
formulas (see
\cite{Luz8})
\begin{equation}\label{koefv_N_diskr}
v_{\mu,N}(k)=\sum_{l=[-\frac{k}{\mu}]'}^{\min\{
[\frac{N-k}{\mu}],n\}}(-1)^l
{n \choose l}b_{\mu, N}(l\mu+k),\quad k=-\mu n,-\mu n+1,\ldots,-1,
\end{equation}
\begin{equation}\label{koef_N b_diskr}
b_{\mu,N}(k)=\sum_{m=k}^Na(m)d_{\mu}(m-k)=\bigl(D^{\mu}_{N}\me
a_N\bigr)_k,\quad
k=0,1,\ldots,N,
\end{equation}
where by $[x]'$ we denote the least integer number among the numbers
that are greater than or equal to $x$, the coefficients $\{d_{\mu
}(k):k\geq0\}$ are determined by the relationship
\[
\sum_{k=0}^{\infty}d_{\mu}(k)x^k=\Bigg(\sum_{j=0}^{\infty}x^{\mu
j}\Bigg)^n,
\]
the matrix $D^{\mu}_{N}$ of dimension $(N+1)\times(N+1)$ is defined by
the coefficients $(D^{\mu}_{N})_{k,j}=d_{\mu}(j-k)$ if
$0\leq k\leq j\leq N$, and $(D^{\mu}_{N})_{k,j}=0$ if $0\leq j<k\leq
N$; and $\me a_N=(a(0),a(1),a(2),\ldots,a(N))'$ is a vector of
dimension $(N+1)$.

The functional $H_N\xi$ from representation (\ref{zobrazh
A_N_i_st.n_d}) has finite variance, and the functional $V_N\zeta$
depends on the known observations of the stochastic sequence $\zeta(k)$
at the points $k=-\mu n,-\mu n+1,\ldots,-1$. Therefore, optimal estimates
$\widehat{A}_N\xi$ and $\widehat{H}_N\xi$ of the functionals ${A}_N\xi$ and ${H}_N\xi$  and the mean-square errors
$\varDelta(f,g;\widehat{A}_N\xi)=\mt E |A_N\xi-\widehat{A}_N\xi
|^2$ and
$\varDelta(f,g;\widehat{H}_N\xi)=\mt E
|H_N\xi-\widehat{H}_N\xi|^2$ of the estimates $\widehat{A}_N\xi$ and
$\widehat{H}_N\xi$ satisfy the following relations:
\begin{align}
\label{mainformula_i_st.n_d}
\widehat{A}_N\xi&=\widehat{H}_N\xi-V_N\zeta,\notag\\
\varDelta(f,g;\widehat{A}_N\xi)&=\mt E |A_N\xi-\widehat{A}_N\xi|^2=\mt E|H_N\xi-V_N\zeta-\widehat{H}_N\xi+V_N\zeta|^2\notag\\
&=\mt E|H_N\xi-\widehat{H}_N\xi|^2=\varDelta(f,g;\widehat{H}_N\xi).
\end{align}
Thus, the interpolation problem for the functional $A_N\xi$ is
equivalent to the interpolation problem for the functional $H_N\xi$.
This problem can be solved by applying the Hilbert space projection
method proposed by Kolmogorov \cite{Kolmogorov}. The optimal linear
estimate $\widehat{A}_N\xi$ of the functional $A_N\xi$ can be
represented in the form
\begin{equation}\label{otsinka A_i_st.n_d}
\widehat{A}_N\xi=\ip
h_{\mu}(\lambda)dZ_{\xi^{(n)}+\eta^{(n)}}(\lambda)-\sum_{k=-\mu
n}^{-1}v_{\mu,N}(k)\big(\xi(k)+\eta(k)\big),
\end{equation}
where
$h_{\mu}(\lambda)$ is the spectral characteristic of the optimal estimate
$\widehat{H}_N\xi$.

Let $H^{0-}(\xi^{(n)}_{\mu}+\eta^{(n)}_{\mu})$ be the closed linear
subspace generated by elements\break$\{\xi^{(n)}(k,\mu)+\eta
^{(n)}(k,\mu
):k\leq-1\}$
of the Hilbert space $H=L_2(\varOmega,\mathcal{F},\mt P)$ of random
variables $\gamma$ with zero mean value and finite variance, $\mt
E\gamma=0$, $\mt E|\gamma|^2<\infty$, with the inner product
$(\gamma
_1;\gamma_2)=\mt E\gamma_1\overline{\gamma_2}$.
Let $H^{N+}(\xi^{(n)}_{-\mu}+\eta^{(n)}_{-\mu})$ be the~closed linear
subspace of the Hilbert space $H=L_2(\varOmega,\mathcal{F},\mt P)$
generated by elements
$\{\xi^{(n)}(k,-\mu)+\eta^{(n)}(k,-\mu):k\geq N+1\}$. The equality
$\xi
^{(n)}(k,-\mu)=(-1)^n\xi^{(n)}(k+\mu n,\mu)$ implies
\[
H^{N+}\bigl(\xi_{-\mu}^{(n)}+\eta_{-\mu}^{(n)}\bigr)=H^{(N+\mu
n)+}\bigl(\xi_{\mu
}^{(n)}+\eta_{\mu}^{(n)}\bigr).
\]
Let us also define the subspaces $L_2^{0-}(p)$ and $L_2^{N+}(p)$ of the
Hilbert space $L_2(p)$ with the inner product $(x_1;x_2)=\ip
x_1(\lambda
)\overline{x_2(\lambda)}p(\lambda)d\lambda$ that are generated by
the functions
$\{e^{i\lambda
k}(1-e^{-i\lambda\mu})^n(i\lambda)^{-n}:k\leq-1\}$ and
$\{e^{i\lambda
k}(1-e^{-i\lambda\mu})^n (i\lambda)^{-n}:k\geq N+1\}$, respectively,
where the function
\[p(\lambda)=f(\lambda)+\lambda^{2n}g(\lambda)\]
is
the spectral density of the sequence $\zeta(m)$, $m\in\mr Z$ \cite{Luz8}.
The optimal estimate
$\widehat{H}_N\xi$ of the functional $H_N\xi$ is the projection of the
element $H_N\xi$ of the Hilbert space $H=L_2(\varOmega,\mathcal
{F},\mt P)$
onto the subspace
\[
H^{0-}\bigl(\xi^{(n)}_{\mu}+\eta^{(n)}_{\mu}\bigr)\oplus
H^{N+}\bigl(\xi^{(n)}_{-\mu}+\eta^{(n)}_{-\mu}\bigr)=H^{0-}\bigl
(\xi^{(n)}_{\mu}+\eta
^{(n)}_{\mu}\bigr)\oplus H^{(N+\mu n
)+}\bigl(\xi^{(n)}_{\mu}+\eta^{(n)}_{\mu}\bigr).
\]
The following conditions characterize the estimate $\widehat{H}_N\xi$:
\begin{enumerate}
\item[1)]
$ \widehat{H}_N\xi\in H^{0-}(\xi^{(n)}_{\mu}+\eta^{(n)}_{\mu})\oplus H^{(N+\mu n)+}(\xi^{(n)}_{\mu}+\eta^{(n)}_{\mu}) $;
\item[2)] $(H_N\xi-\widehat{H}_N\xi)\perp H^{0-}(\xi^{(n)}_{\mu}+\eta^{(n)}_{\mu})\oplus H^{(N+\mu n)+}(\xi^{(n)}_{\mu}+\eta^{(n)}_{\mu})$.
\end{enumerate}
The functional $H_N\xi$ in the space $H$ admits the spectral representation
\[
H_N\xi=\int_{-\pi}^{\pi}B^{\mu}_N\bigl(e^{i\lambda}\bigr)\frac
{(1-e^{-i\lambda\mu})^n}{(i\lambda)^n}
dZ_{\xi^{(n)}+\eta^{(n)}}(\lambda)
-\int_{-\pi}^{\pi}A_N\bigl(e^{i\lambda}\bigr)dZ_{\eta}(\lambda),
\]
\[
B^{\mu}_N\bigl(e^{i\lambda}\bigr)=\sum_{k=0}^{N}b_{\mu
,N}(k)e^{i\lambda k}
=\sum_{k=0}^{N}\bigl(D^{\mu}_N\me a_N\bigr)_ke^{i\lambda k},\qquad
A_N\bigl(e^{i\lambda}\bigr)
=\sum_{k=0}^{N}a(k)e^{i\lambda k}.
\]
Making use of the described representation and condition 2), we derive
the following equation for determining the spectral characteristic
$h_{\mu}(\lambda)$:
\begin{align*}
&\ip\bigg[\bigg(B^{\mu}_N\bigl(e^{i\lambda}\bigr)\frac{(1-e^{-i\lambda\mu})^n}{(i\lambda)^n}-h_{\mu}(\lambda)\bigg)p(\lambda)-A\bigl(e^{i\lambda}\bigr)g(\lambda)(-i\lambda)^n\bigg]\\
&\quad \times\frac{(1-e^{i\lambda\mu})^n}{(-i\lambda)^{n}}e^{-i\lambda k}d\lambda=0\quad\forall k\leq-1,\,\forall k\geq N+\mu n+1.
\end{align*}
Thus, the spectral characteristic
$h_{\mu}(\lambda)$ can be represented as follows:
\begin{eqnarray}\label{spectr A_i_st.n_d}
&\displaystyle h_{\mu}(\lambda)=B^{\mu}_N\bigl(e^{i\lambda
}\bigr)\frac{(1-e^{-i\lambda\mu})^n}{(i\lambda
)^n}-A_N\bigl(e^{i\lambda}\bigr)
\frac{(-i\lambda)^ng(\lambda)}{p(\lambda)}-\frac{(-i\lambda)^{n}
C^{\mu}_N(e^{i\lambda})}{(1-e^{i\lambda\mu})^np(\lambda)},&\notag\\
&\displaystyle C^{\mu}_N\bigl(e^{i\lambda}\bigr)=\sum_{k=0}^{N+\mu n}c_{\mu}(k)e^{i\lambda k},&
\end{eqnarray}
where $c_{\mu}(k)$, $k=0,1,2,\ldots,N+\mu n$, are unknown coefficients
we have to determine. Condition 1) implies that the spectral
characteristic $h_{\mu}(\lambda)$ satisfies the following equations:
\begin{align*}
&\ip\ld[B^{\mu}_N\bigl(e^{i\lambda}\bigr)-
\frac{A_N(e^{i\lambda})\lambda^{2n}g(\lambda)}
{(1-e^{-i\lambda\mu})^np(\lambda)}-
\frac{\lambda^{2n}C^{\mu}_N(e^{i\lambda})}
{|1-e^{i\lambda\mu}|^{2n}p(\lambda)}
\rd]e^{-i\lambda l}d\lambda=0,\\
&\quad  0\leq l\leq N+\mu n.
\end{align*}
The derived equations are represented as a system of $N+\mu n+1$ linear
equations:
\begin{equation}\label{linear equations1}
b_{\mu,N}(l)-\sum_{ m=0}^{N+\mu n}T^{\mu}_{l,m}a_{\mu,N}(m)
=\sum_{k=0}^{N+\mu n}P_{l,k}^{\mu}c_{\mu}(k),\quad0\leq l\leq N,
\end{equation}
\begin{equation}\label{linear equations2}
-\sum_{ m=0}^{N+\mu n}T^{\mu}_{l,m}a_{\mu,N}(m)
=\sum_{k=0}^{N+\mu n}P_{l,k}^{\mu}c_{\mu}(k),\quad N+1\leq l\leq
N+\mu n,
\end{equation}
where the coefficients $\{a_{\mu,N}(m):0\leq m\leq N+\mu n\}$ are
calculated by the formula
\begin{equation}\label{coeff a_N_mu}
a_{\mu,N}(m)
=\sum_{l=\max\{[\frac{m-N}{\mu}]',0\}}^{\min\{[\frac
{m}{\mu}],n\}}
(-1)^l{n \choose l}a(m-\mu l),\quad0\leq m\leq N+\mu n,
\end{equation}
and the Fourier coefficients $\{T^{\mu}_{k,j},P_{k,j}^{\mu}:0\leq k,j
\leq N+\mu n\}$ are calculated by the formulas
\[
T^{\mu}_{k,j}=\frac{1}{2\pi}\ip
e^{i\lambda(j-k)}
\frac{\lambda^{2n}g(\lambda)}
{|1-e^{i\lambda\mu}|^{2n}(f(\lambda)+\lambda^{2n}g(\lambda))}
d\lambda,\quad
0\leq k,j \leq
N+\mu n,
\]
\[
P_{k,j}^{\mu}=\frac{1}{2\pi}\ip e^{i\lambda(j-k)}
\dfrac{\lambda^{2n}}
{|1-e^{i\lambda\mu}|^{2n}(f(\lambda)+\lambda^{2n}g(\lambda))}
d\lambda,\quad0\leq k,j \leq N+\mu n.
\]

Denote by $[D_N^{\mu}\me a_N]_{+\mu n}$ the vector of dimension
$(N+\mu
n +1)$ constructed by adding
$\mu n$ zeros to the vector $D_N^{\mu}\me a_N$ of dimension $(N+1)$.
Using these definitions, system (\ref{linear equations1})--(\ref{linear
equations2}) can be represented in the matrix form
\begingroup
\abovedisplayskip=9pt
\belowdisplayskip=9pt
\[
\bigl[D_N^{\mu}\me a_N\bigr]_{+\mu n}-\me T^{\mu}_N\me a^{\mu
}_N=\me P^{\mu}_N\me
c^{\mu}_N,
\]
where
\[
\me a^{\mu}_N=\big(a_{\mu,N}(0),a_{\mu,N}(1),a_{\mu,N}(2),\ldots,a_{\mu ,N}(N+\mu n)\big)'
\]
and
\[
\me c^{\mu}_N=\big(c_{\mu}(0),c_{\mu}(1),c_{\mu}(2), \ldots,c_{\mu}(N+\mu n)\big)'
\]
are vectors of dimension $(N+\mu n+1)$; and $\me P^{\mu}_N$ and $\me
T^{\mu}_N$ are matrices of dimension $(N+\mu n+1)\times(N+\mu n+1)$
with elements $(\me P^{\mu}_N)_{l,k}=P_{l,k}^{\mu}$ and $(\me T^{\mu
}_N)_{l, k} =T^{\mu}_{l,k}$, $0\leq l,k\leq N+\mu n$.
Thus, the coefficients $c_{\mu}(k)$, $0\leq k\leq N+\mu n$, are
determined by the formula
\[
c_{\mu}(k)=\bigl(\bigl(\me P^{\mu}_N\bigr)^{-1}\bigl[D_N^{\mu
}\me a_N\bigr]_{+\mu n}-\bigl(\me
P^{\mu}_N\bigr)^{-1}\me T^{\mu}_N\me a_{\mu}\bigr)_k,\quad0\leq
k\leq N+\mu n,
\]
where $((\me P^{\mu}_N)^{-1}[D_N^{\mu}\me a_N]_{+\mu n}-(\me P^{\mu
}_N)^{-1}\me T^{\mu}_N\me a^{\mu}_N)_k$, $0\leq k\leq N+\mu n$, is the
$k$th element of the vector $(\me P^{\mu}_N)^{-1}[D_N^{\mu}\me
a_N]_{+\mu n}-(\me P^{\mu}_N)^{-1}\me T^{\mu}_N\me a^{\mu}_N$.
The existence of the invertible matrix $(\me P^{\mu}_N)^{-1}$ was shown
in \cite{Salehi} under condition (\ref{umova11_i_st.n_d}). The spectral
characteristic $h_{\mu}(\lambda)$ of the estimate $\widehat{H}_N\xi
$ of
the functional $H_N\xi$ is calculated by formula (\ref{spectr
A_i_st.n_d}), where
\[
C^{\mu}_N\bigl(e^{i\lambda}\bigr)=\sum_{k=0}^{N+\mu n}
\bigl(\bigl(\me P^{\mu}_N\bigr)^{-1}\bigl[D_N^{\mu}\me a_N\bigr
]_{+\mu n}-\bigl(\me P^{\mu
}_N\bigr)^{-1}\me T^{\mu}_N\me a^{\mu}_N\bigr)_k e^{i\lambda
k}.
\]
The value of the mean-square errors of the estimates $\widehat{A}_N\xi$ and $\widehat
{H}_N\xi$ can be calculated by the formula
\begin{align}
\varDelta(f,g;\widehat{A}_N\xi)&=\varDelta(f,g;\widehat{H}_N\xi)=\mt E|H_N\xi-\widehat{H}_N\xi|^2\notag\\[4pt]
&=\frac{ 1}{2\pi}\ip\frac{|A_N(e^{i\lambda})(1-e^{i\lambda\mu})^{n}f(\lambda)-\lambda^{2n}C^{\mu}_N(e^{i\lambda})|^2}{|1-e^{i\lambda\mu}|^{2n}(f(\lambda)+\lambda^{2n}g(\lambda))^2}g(\lambda)d\lambda\notag\\[4pt]
&\quad +\frac{ 1}{2\pi}\ip\frac{|A_N(e^{i\lambda})(1-e^{i\lambda\mu})^{n}\lambda^{2n}g(\lambda)+\lambda^{2n}C^{\mu}_N(e^{i\lambda})|^2}{\lambda^{2n}|1-e^{i\lambda\mu}|^{2n}(f(\lambda)+\lambda^{2n}g(\lambda))^2}f(\lambda)d\lambda\notag\\[4pt]
&=\big\langle\bigl[D_N^{\mu}\me a_N\bigr]_{+\mu n}- \me T^{\mu}_N\me a^{\mu}_N,\bigl(\me P^{\mu}_N\bigr)^{-1}\bigl[D_N^{\mu}\me a_N\bigr]_{+\mu n}-\bigl(\me P^{\mu}_N\bigr)^{-1}\me T^{\mu}_N\me a^{\mu}_N\big\rangle\notag\\[4pt]
&\quad +\langle\me Q_N\me a_N,\me a_N\rangle,\label{poh A_i_st.n_d}
\end{align}
\endgroup
where $\me Q_N$ is the matrix of dimension $(N+1)\times(N+1)$ with the
coefficients $(\me Q_N)_{l,k}=Q_{l,k}$, $0\leq l,k\leq N$, calculated
by the formula
\[
Q_{k,j}=\frac{1}{2\pi}\ip
e^{i\lambda(j-k)}\frac{f(\lambda)g(\lambda)}
{f(\lambda)+\lambda^{2n}g(\lambda)}d\lambda,\quad0\leq k,j \leq N.
\]

We can summarize the derived results in the form of the following theorem.

\begin{thm}\label{thm1_i_st.n_d}
Let $\{\xi(m),m\in\mr Z\}$ be a stochastic sequence with stationary
$n$th increments $\xi^{(n)}(m,\mu)$, and let $\{\eta(m),m\in\mr Z\}
$ be
a stationary stochastic sequence uncorrelated with $\xi(m)$. Let the
spectral densities $f(\lambda)$ and $g(\lambda)$ of the sequences
satisfy the minimality condition (\ref{umova11_i_st.n_d}). The optimal
linear estimate $\widehat{A}_N\xi$ of the functional $A_N\xi$, which
depends on the values $\xi(m)$, $0\leq m\leq N$, based on the
observations of the sequence
$\xi(m )+\eta(m )$ at points of the set $Z\setminus\{0,1,2,\ldots
,N\}
$ is calculated by formula (\ref{otsinka A_i_st.n_d}). The spectral
characteristic
$h_{\mu}(\lambda)$ and the value of the mean-square error $\varDelta
(f,g;\widehat{A}_N\xi)$ of the optimal estimate $\widehat{A}_N\xi$ are
calculated by formulas (\ref{spectr A_i_st.n_d}) and (\ref{poh
A_i_st.n_d}), respectively.
\end{thm}

\begin{cor}\label{nas_A_xi_i_st.n_d}
Let the spectral density $f(\lambda)$ of the sequence $\xi(m)$ satisfy
the minimality condition
\[
\ip\frac{\lambda^{2n}}{|1-e^{i\lambda\mu}|^{2n}f(\lambda)}
d\lambda<\infty.
\]
The optimal linear estimate $\widehat{A}_N\xi$ of the functional
$A_N\xi
$ of unknown values $\xi(m)$, $0\leq m\leq N$, based on observations of
the sequence
$\xi(m)$ at the points $m\in\mr Z\setminus\{0,1,2,\ldots,N\}$ can be
calculated by the formula
\begin{equation}\label{estim A_xi_i_st.n_d}
\widehat{A}_N\xi=\ip
h_{\mu}^{\xi}(\lambda)dZ_{\xi^{(n)}}(\lambda)-\sum_{k=-\mu
n}^{-1}v_{\mu,N}(k)\xi(k).
\end{equation}
The spectral characteristic
$h_{\mu}^{\xi}(\lambda)$ and the mean-square error $\varDelta
(f;\widehat
{A}_N\xi)$ of the optimal estimate $\widehat{A}_N\xi$ can be calculated
by the formulas
\begin{equation}\label{spectr_A_xi_i_st.n_d}
h_{\mu}^{\xi}(\lambda)=B^{\mu}_N\bigl(e^{i\lambda}\bigr)
\frac{(1-e^{-i\lambda\mu})^n}
{(i\lambda)^n}
-\frac{(-i\lambda)^{n}\sum_{k=0}^{N+\mu n}((\me F^{\mu}_N)^{-1}[D_N^{\mu}\me a_N]_{+\mu n})_k e^{i\lambda k}}
{(1-e^{i\lambda\mu})^nf(\lambda)},
\end{equation}
\begin{align}
\varDelta(f;\widehat{A}_N\xi)&=\frac{1}{2\pi}\ip
\frac{\lambda^{2n}|\sum_{k=0}^{N+\mu n} ((\me F^{\mu}_N)^{-1}[D_N^{\mu}\me a_N]_{+\mu n})_ke^{i\lambda k}|^2}
{|1-e^{i\lambda\mu}|^{2n}f(\lambda)}
d\lambda\notag
\\
&=\big\langle\bigl(\me F^{\mu}_N\bigr)^{-1}\bigl[D_N^{\mu}\me
a_N\bigr]_{+\mu n},
\bigl[D_N^{\mu}\me a_N\bigr]_{+\mu n}\big\rangle,\label{poh A_xi_i_st.n_d}
\end{align}
where $\me F^{\mu}_N$ is the matrix of dimension $(N+\mu n+1)\times
(N+\mu n+1)$ with elements $(\me F^{\mu}_N)_{k,j}=F^{\mu}_{k,j}$,
$0\leq k,j\leq N+\mu n$,
\[
F_{k,j}^{\mu}=\frac{1}{2\pi}\ip e^{i\lambda(j-k)}
\dfrac{\lambda^{2n}}{|1-e^{i\lambda\mu}|^{2n}f(\lambda)}
d\lambda,\quad0\leq k,j\leq N+\mu n.
\]
\end{cor}

In the case of estimation of an unobserved value $\xi(p)$, $0\leq
p\leq
N$, the following statement holds true.

\begin{thm}\label{thm xi_i_st.n_d}
Let the conditions of Theorem \ref{thm1_i_st.n_d} hold. The optimal
linear estimate $\widehat{\xi}(p)$ of an unobserved value
$\xi(p)$, $0\leq p\leq N$, of the stochastic sequence with $n$th
stationary increments based on observations of the sequence $\xi
(m)+\eta
(m)$ at the points $m\in\mr Z\setminus\{0,1,2,\ldots,N\}$ is calculated
by the formula
\[
\widehat{\xi}(p)=\ip
h_{\mu,p}(\lambda)dZ_{\xi^{(n)}+\eta^{(n)}}(\lambda)-\sum
_{l=1}^{n}(-1)^l{n \choose l}\bigl(\xi(p-\mu l)+\eta(p-\mu l)\bigr),
\]
\[
h_{\mu,p}(\lambda)=\frac{(1-e^{-i\lambda\mu}
)^n}
{(i\lambda)^n}\sum
_{k=0}^pd_{\mu}(p-k)e^{i\lambda k}-
\frac{e^{i\lambda p }(-i\lambda)^ng(\lambda)}{p(\lambda)}
-\frac{(-i\lambda)^{n}C_p^{\mu}(e^{i\lambda})}
{(1-e^{i\lambda\mu}
)^np(\lambda)},
\]
\[
C_p^{\mu}(e^{i\lambda})=\sum_{k=0}^{N+\mu n}
\bigl(\bigl(\me P^{\mu}_N\bigr)^{-1}\me d_{\mu,p}-\bigl(\me
P^{\mu}_N\bigr)^{-1}\me T^{\mu
}_p\me a_{n}\bigr)_k e^{i\lambda
k},
\]
where
\[
\me d_{\mu,p}=\bigl(d_{\mu}(p),d_{\mu}(p-1),d_{\mu}(p-2),\ldots
,d_{\mu
}(0),0,\ldots,0\bigr)'
\]
and
\[
\me a_{n}=\bigl(a_{n}(0),a_{n}(1), \ldots,a_{n}(n),0,\ldots,0\bigr)',
\]
$a_{n}(k)=(-1)^k{n \choose k}$, $k=0,1,2,\ldots,n$, are vectors of
dimension $(N+\mu n+1)$, $\me T^{\mu}_p$ is the $(N+\mu n +1)\times
(N+\mu n+1)$ matrix with elements
$(\me T^{\mu}_p)_{l, k} =T^{\mu}_{l,p+\mu k}$ if $0\leq l\leq N+\mu
n$, $0\leq k\leq n$, and $(\me T^{\mu}_p)_{l, k} =0$ if $0\leq l\leq
N+\mu n$, $N+1\leq k\leq N+\mu n$.
The value of the mean-square error of the optimal estimate is
calculated the by formula
\[
\varDelta\big(f,g;\widehat{\xi}(p)\big)= \big\langle\me d_{\mu,p}-\me
T^{\mu}_p \me
a_{n},\bigl(\me P^{\mu}_N\bigr)^{-1}\me d_{\mu,p}-\bigl(\me P^{\mu
}_N\bigr)^{-1}\me T^{\mu
}_p \me a_{n}\big\rangle+ Q_{0,0}.
\]
\end{thm}

\begin{cor}\label{nas xi_i_st.n_d}
In the case of estimating the sequence $\xi(m)$ with $n$th stationary
increments at points of the set $\mr Z\setminus\{0,1,2,\ldots,N\}$, the
optimal linear estimate of a value
$\xi(p)$, $0\leq p\leq N$, is calculated by the formula
\[
\widehat{\xi}(p)=\ip
h^{\xi}_{\mu,p}(\lambda)dZ_{\xi^{(n)}}(\lambda)-\sum_{l=1}^{n}(-1)^l{n
\choose l}\xi(p-\mu l),
\]
\[
h^{\xi}_{\mu,p}(\lambda)=
\frac{(1-e^{-i\lambda\mu})^n}
{(i\lambda)^n}\sum
_{k=0}^pd_{\mu}(p-k)e^{i\lambda k}-
\frac{(-i\lambda)^{n}\sum_{k=0}^{N+\mu n}((\me F^{\mu}_N)^{-1}\me d_{\mu,p})_k e^{i\lambda k}}
{(1-e^{i\lambda\mu})^nf(\lambda)}.
\]
The value of the mean-square error of the estimate is calculated by the formula
\begin{align*}
\varDelta\big(f;\widehat{\xi}(p)\big)&=
\frac{ 1}{2\pi}\ip
\frac{\lambda^{2n}|\sum_{k=0}^{N+\mu n} ((\me F^{\mu}_N)^{-1}\me d_{\mu,p})_ke^{i\lambda k}|^2}
{|1-e^{i\lambda\mu}|^{2n}f^2(\lambda)}
\,d\lambda\notag
\\
&= \big\langle\bigl(\me F^{\mu}_N\bigr)^{-1}\me d_{\mu,p},\me
d_{\mu,p}\big\rangle.
\end{align*}
\end{cor}

\begin{exm}
Consider the stochastic sequence $\xi(m)$, $m\in\mr Z$, defined by
the equation
\[
\xi(m)=(1-\phi)\xi(m-1)+\phi\xi(m-2)+\varepsilon(m),
\]
which means that values of the sequence $\xi(m)$ are defined as a
weighted sum of two previous values of the sequence plus a value
$\varepsilon(m)$ of the sequence of independent identically distributed
random variables with mean value $\mt E\varepsilon(m)=0$ and variance
$\mt E \varepsilon^2(m)=1$.

Consider the increment $\xi^{(1)}(m;1)=\xi(m)-\xi(m-1)$ of the
sequence. We can find that
\[
\xi^{(1)}(m;1)=-\phi\xi^{(1)}(m-1;1)+\varepsilon(m).
\]
Thus, the increment sequence $\xi^{(1)}(m;1)$ with step $\mu=1$ is an
autoregressive sequence with parameter $0<\phi<1$. The sequence $\xi
(m)$ is an ARIMA(1;1;0) sequence with the spectral density
\[
f(\lambda)=\frac{\lambda^2}
{|1-e^{-i\lambda}|^2|1+\phi e^{-i\lambda}|^2}.
\]

Let us find the estimate $\widehat{A}_1\xi$ of the value of the
functional $A_1\xi=2\xi(0)+\xi(1)$ based on observations of the
sequence $\xi(m)$ at the points $m\in\mr Z\setminus\{0,1\}$.
Let $\phi={1}/{2}$. In this case, $v_{1,1}(-1)=-2$,
\[
\me F_1\,{=}\,\frac{1}{4}
\begin{pmatrix}
5&2&0\\
2&5&2\\
0&2&5\\
\end{pmatrix}
,\qquad
\me F^{-1}_1\,{=}\,\frac{4}{85}
\begin{pmatrix}
21&-10&4\\
-10&25&-10\\
4&-10&21\\
\end{pmatrix}
,\qquad
\bigl[D_1^{1}\me a_1\bigr]_{+1}\,{=}
\begin{pmatrix}
3\\
1\\
0\\
\end{pmatrix}
\!\!.
\]
Therefore,
\[
h_1^{\xi}(\lambda)=-\frac{106}{85}e^{-i\lambda}-\frac
{4}{85}e^{3i\lambda},
\]
\begin{align*}
\widehat{A}_1\xi&=-\frac{106}{85}\xi^{(1)}(-1;1)-\frac{4}{85}\xi
^{(1)}(3;1)-3\xi(-1)\\
&=\frac{106}{85}\xi(-2)+\frac{149}{85}\xi(-1)+\frac{4}{85}\xi
(2)-\frac
{4}{85}\xi(3).
\end{align*}
The value of the mean-square error of the estimate is $\varDelta
(f,g;\widehat{A}_1\xi)=\tfrac{88}{17}$.
\end{exm}

\section{Interpolation of cointegrated sequences}

Consider two integrated sequences $\{\xi(m),m\in\mr Z\}$ and $\{\zeta
(m),m\in\mr Z\}$ with absolutely continuous spectral functions
$F(\lambda)$ and $P(\lambda)$ and the corresponding spectral densities
$f(\lambda)$ and $p(\lambda)$.

\begin{defin}\label{ozn_coint}
Two integrated sequences $\{(\xi(m),\zeta(m)),m\in\mr Z\}$ are called
cointegrated (of order $0$) if, for some constant $\beta\neq0$, the
linear combination
$\{\zeta(m)-\beta\xi(m): m\in\mr Z\}$ is a stationary sequence.
\end{defin}

The interpolation problem for cointegrated sequences consists in
mean-square optimal linear estimation of the functional
\[
A_N\xi=\sum_{k=0}^{N}a(k)\xi(k)
\]
of unknown values of the stochastic sequence $\xi(m)$ based on
observations of the stochastic sequence $\zeta(m)$ at the points $m\in
\mr Z\setminus\{0,1,2,\ldots,N\}$. To solve the problem, we can use the
results obtained in the previous sections.

Suppose that the spectral density
$p(\lambda)$ of the sequence $\zeta(m)$ satisfies the minimality
condition
\begin{equation}
\ip\frac{\lambda^{2n}}
{|1-e^{i\lambda\mu}|^{2n}p(\lambda)}
d\lambda<\infty.\label{umova111_i_st.n_d}
\end{equation}
Let the matrices $\me P^{\mu,\beta}_N$, $\me T_{N}^{\mu,\beta}$,
$\me
Q_{N}^{\beta}$ be defined by the Fourier coefficients of the functions
\begin{equation}\label{functions for
beta-operators}
\frac{\lambda^{2n}}
{|1-e^{i\lambda\mu}|^{2n}p(\lambda)},
\qquad
\dfrac{p(\lambda)-\beta^2f(\lambda)}{|1-e^{i\lambda\mu}
|^{2n}p(\lambda
)},\qquad
\dfrac{[f(\lambda)p(\lambda)-\beta^2f^2(\lambda)
]_+}{\lambda^{2n}
p(\lambda)}
\end{equation}
in the same way as the matrices $\me P^{\mu}_N$, $\me T^{\mu}_{N}$,
$\me Q_{N}$ were defined.
Theorem \ref{thm1_i_st.n_d} implies the following formula for
calculating the spectral characteristic $h^{\beta}_{\mu,N}(\lambda)$
of the optimal estimate
\begin{equation}\label{otsinka A_N_coi_i_st.n_d} \widehat{A}_N\xi
=\ip
h^{\beta}_{\mu,N}(\lambda)dZ_{\zeta^{(n)}}(\lambda)-\sum_{k=-\mu
n}^{-1}v_{\mu,N}(k)\zeta(k)
\end{equation}
of the functional $A_N\xi$:
\begin{equation}\label{spectr A_N_coi_i_st.n_d}
h^{\beta}_{\mu,N}(\lambda)=B_N^{\mu}\big(e^{i\lambda}\big)
\frac{(1-e^{-i\lambda \mu})^n}
{(i\lambda)^n}-
A_N\bigl(e^{i\lambda}\bigr)
\frac{p(\lambda)
- \beta^2f(\lambda)}
{(i\lambda)^np(\lambda)}
-\frac{(-i\lambda)^{n}C_{\mu,N}^{\beta}(e^{i\lambda})}{(1-e^{i\lambda
\mu})^np(\lambda)},
\end{equation}
where
\[
C_{\mu,N}^{\beta}\bigl(e^{i\lambda}\bigr)=\sum_{k=0}^{N+\mu
n}\big(\bigl(\me P^{\mu
,\beta}_{N}\bigr)^{-1}\bigl[D_N^{\mu}\me a_N\bigr]_{+\mu n}-
\bigl(\me P^{\mu,\beta}_{N}\bigr)^{-1}\me T^{\mu,\beta}_{N}\me
a^{\mu}_{N}\big)_ke^{i\lambda
k}.
\]

The value of the mean-square error of the estimate $\widehat{A}_N\xi$
is calculated by the formula
\begin{align}
&\varDelta(f,g;\widehat{A}_N\xi)\notag\\
&\quad =\frac{ 1}{2\pi}\ip\frac{|A_N(e^{i\lambda})(1-e^{i\lambda\mu})^{n}\beta^2f(\lambda)-\lambda^{2n}C_{\mu,N}^{\beta}(e^{i\lambda})|^2}{\lambda^{2n}|1-e^{i\lambda\mu}|^{2n}p^2(\lambda)}p(\lambda)d\lambda\notag\\
&\qquad -\frac{ \beta^2}{2\pi}\ip\frac{|A_N(e^{i\lambda})(1-e^{i\lambda\mu})^{n}\beta^2f(\lambda)-\lambda^{2n}C_{\mu,N}^{\beta}(e^{i\lambda})|^2}{\lambda^{2n}|1-e^{i\lambda\mu}|^{2n}p^2(\lambda)}f(\lambda)d\lambda\notag\\
&\qquad +\frac{ 1}{2\pi}\ip\frac{|A_N(e^{i\lambda})(1-e^{i\lambda\mu})^{n}[p(\lambda)-\beta^2f(\lambda)]_++\lambda^{2n}C_{\mu,N}^{\beta}(e^{i\lambda})|^2}{\lambda^{2n}|1-e^{i\lambda\mu}|^{2n}p^2(\lambda)}f(\lambda)d\lambda\notag\\
&\quad = \big\langle\bigl[D_N^{\mu}\me a_N\bigr]_{+\mu n}-\me T^{\mu,\beta}_{N}\me a^{\mu}_{N},\bigl(\me P^{\mu,\beta}_{N}\bigr)^{-1}\bigl[D_N^{\mu}\me a_N\bigr]_{+\mu n}-\bigl(\me P^{\mu,\beta}_{N}\bigr)^{-1}\me T^{\mu,\beta}_{N}\me a^{\mu}_{N}\big\rangle\notag\\
&\qquad +\big\langle\me Q^{\beta}_N\me a_N,\me a_N \big\rangle.\label{poh A_N_coi_i_st.n_d}
\end{align}

The described results are presented as the following theorem.

\begin{thm}\label{thm5_i_st.n_d}
Let $\{(\xi(m),\zeta(m)),m\in\mr Z\}$ be two cointegrated sequences
with spectral densities
$f(\lambda)$ and $p(\lambda)$, and let the spectral density
$p(\lambda)$ satisfy the minimality condition
\eqref{umova111_i_st.n_d}. If the stochastic sequences $\xi(m)$ and
$\zeta(m)-\beta\xi(m)$ are uncorrelated, then the spectral
characteristic $h^{\beta}_{\mu,N}(\lambda)$ and the value of the
mean-square error $\varDelta(f,g;\widehat{A}_N\xi)$ of the optimal
estimate $\widehat{A}_N\xi$ \eqref{otsinka A_N_coi_i_st.n_d} of the
functional $A_N\xi$ based on the observations of the sequence
$\zeta(m)$ at the points $m\in\mr Z\setminus\{0,1,2,\ldots,N\}$ are
calculated by formulas \eqref{spectr A_N_coi_i_st.n_d} and
\eqref{poh A_N_coi_i_st.n_d}, respectively.
\end{thm}

\section{Minimax-robust method of interpolation}

Formulas for calculating values of the mean-square error $\varDelta
(h(f,g);f,g)=\break\varDelta(f,g;\widehat{A}_N\xi)=\mt E |A_N\xi
-\widehat{A}_N\xi
|^2$ and the spectral characteristics of the optimal estimates
of the functional ${A}_N\xi$ based on observations of the sequence
$\xi(m)+\eta(m)$ can be applied under the condition that
the spectral densities
$f(\lambda)$ and $g(\lambda)$ of the stochastic sequences $\xi(m)$ and
$\eta(m)$ are known.
However, these formulas often cannot be used in many practical
situations since the exact values of the densities are not available.
In this situation, the minimax-robust method can be applied. It
consists in finding the estimate that provides
a~minimum of the mean-square errors for all spectral densities from a
given set
$\md D=\md D_f\times\md D_g$ of admissible spectral densities simultaneously.

\begin{defin} For a given class of spectral densities $\mathcal{D}=\md
D_f\times\md D_g$, spectral densities
$f^0(\lambda)\in\mathcal{D}_f$ and $g^0(\lambda)\in\md D_g$ are
called the
least favorable densities in the class $\mathcal{D}$ for the optimal linear
interpolation of the functional $A_N\xi$ if the following relation holds:
\begingroup
\abovedisplayskip=7.5pt
\belowdisplayskip=7.5pt
\[
\varDelta(f^0,g^0)=\varDelta\bigl(h(f^0,g^0);f^0,g^0\bigr)=
\max_{(f,g)\in\mathcal{D}_f\times\md
D_g}\varDelta\bigl(h(f,g);f,g\bigr).
\]
\end{defin}

\begin{defin}
For a given class of spectral densities $\mathcal{D}=\md
D_f\times\md D_g$, the spectral characteristic $h^0(\lambda)$ of
the optimal linear estimate of the functional $A_N\xi$ is called
minimax-robust if the following conditions are satisfied:
\[
h^0(\lambda)\in H_{\mathcal{D}}=\bigcap_{(f,g)\in\mathcal
{D}_f\times\md D_g}
L_2^{0- }(p)\oplus L_2^{(N+\mu n)+ }(p),\vspace*{-9pt}
\]
\[
\min_{h\in H_{\mathcal{D}}}\max_{(f,g)\in\mathcal{D}_f\times\md
D_g}\varDelta(h;f,g)=\max_{(f,g)\in\mathcal{D}_f\times\md
D_g}\varDelta\bigl(h^0;f,g\bigr).
\]
\end{defin}

\begin{lemma}
The spectral densities $f^0\in\mathcal{D}_f$ and
$g^0\in\mathcal{D}_g$ that satisfy the minimality condition \eqref
{umova11_i_st.n_d}
are the least favorable in the class $\mathcal{D}$ for the optimal
linear interpolation of the functional $A_N\xi$ based on observations
of the sequence $\xi(m)+\eta(m)$
at the points $m\in\mr Z\setminus\{0,1,2,\ldots,N\}$
if the matrices $ (\me P^{\mu}_N)^0$, $(\me T^{\mu}_N)^0$, $(\me
Q_N)^0$ whose elements are defined by the Fourier coefficients of the functions
\begin{equation}\label{functions_for_lemma_i_st.n_d}
\dfrac{\lambda^{2n}}{|1-e^{i\lambda\mu}|^{2n}p^0(\lambda
)}, \qquad
\dfrac{\lambda^{2n}g^0(\lambda)}{|1-e^{i\lambda\mu}
|^{2n}p^0(\lambda
)},\qquad
\dfrac{f^0(\lambda)g^0(\lambda)}{
p^0(\lambda)},
\end{equation}
where $p^0(\lambda)=f^0(\lambda)+\lambda^{2n}g^0(\lambda)$,
determine a
solution to the constrained optimization problem
\begin{align}
&\max_{(f,g)\in\mathcal{D}_f\times\md D_g}\big(\big\langle\bigl[D_N^{\mu}\me a_N\bigr]_{+\mu n}- \me T^{\mu}_N\me a_{\mu},\bigl(\me P^{\mu}_N\bigr)^{-1}\bigl[D_N^{\mu}\me a_N\bigr]_{+\mu n}-\bigl(\me P^{\mu}_N\bigr)^{-1}\me T^{\mu}_N\me a^{\mu}_N\big\rangle\notag\\
&\qquad +\langle\me Q_N\me a_N,\me a_N\rangle\big)\notag\\
&\quad = \big\langle\bigl[D_N^{\mu}\me a_N\bigr]_{+\mu n}- \bigl(\me T^{\mu}_N\bigr)^0\me a^{\mu}_N,\bigl(\bigl(\me P^{\mu}_N\bigr)^0\bigr)^{-1}\bigl[D_N^{\mu}\me a_N\bigr]_{+\mu n}-\bigl(\bigl(\me P^{\mu}_N\bigr)^0\bigr)^{-1}\bigl(\me T^{\mu}_N\bigr)^0\me a^{\mu}_N\big\rangle\notag\\
&\qquad +\big\langle\me Q^0_N\me a_N,\me a_N\big\rangle.\label{minimax1_i_st.n_d}
\end{align}
The minimax-robust spectral characteristic $h^0=h_{\mu}(f^0,g^0)$ is
calculated by formula \eqref{spectr A_i_st.n_d} if
$h_{\mu}(f^0,g^0)\in H_{\mathcal{D}}$.
\end{lemma}
\endgroup

The presented statements follow from the
introduced definitions and Theorem \ref{thm1_i_st.n_d}.

The minimax-robust spectral characteristic $h^0$ and the least
favorable spectral densities $(f^0,g^0)$
form a saddle point of the function $\varDelta(h;f,g)$ on the set
$H_{\mathcal{D}}\times\mathcal{D}$.
The saddle-point inequalities
\[
\varDelta\bigl(h;f^0,g^0\bigr)\geq\varDelta\bigl(h^0;f^0,g^0\bigr
)\geq\varDelta\bigl(h^0;f,g\bigr)
\quad\forall f\in\mathcal{D}_f,\forall g\in\mathcal{D}_g,\forall
h\in H_{\mathcal{D}}
\]
hold if $h^0=h_{\mu}(f^0,g^0)$,
$h_{\mu}(f^0,g^0)\in H_{\mathcal{D}}$, and $(f^0,g^0)$ is a solution to
the constrained optimization problem
\[
\widetilde{\varDelta}(f,g)=-\varDelta\bigl(h_{\mu}\bigl
(f^0,g^0\bigr);f,g\bigr)\to
\inf,\quad(f,g)\in\mathcal{D},
\]
\begin{align*}
&\varDelta\bigl(h_{\mu}\bigl(f^0,g^0\bigr);f,g\bigr)\\
&\quad =\frac{1}{2\pi}\ip\frac{|A_N(e^{i\lambda})(1-e^{i\lambda\mu})^{n}f^0(\lambda)-\lambda^{2n}C^{\mu,0}_{N}(e^{i\lambda})|^2}{|1-e^{i\lambda\mu}|^{2n}(f^0(\lambda)+\lambda^{2n}g^0(\lambda))^2}g(\lambda)d\lambda\\
&\qquad +\frac{1}{2\pi}\ip\frac{|A_N(e^{i\lambda})(1-e^{i\lambda\mu})^{n}\lambda^{2n}g^0(\lambda)+\lambda^{2n}C^{\mu,0}_{N}(e^{i\lambda})|^2}{\lambda^{2n}|1-e^{i\lambda\mu}|^{2n}(f^0(\lambda)+\lambda^{2n}g^0(\lambda))^2}f(\lambda)d\lambda,
\end{align*}
\[
C^{\mu,0}_{N}\bigl(e^{i\lambda}\bigr)=\sum_{k=0}^{N+\mu n} \bigl
(\bigl(\bigl(\me P^{\mu
}_N\bigr)^0\bigr)^{-1}\bigl[D_N^{\mu}\me a_N\bigr]_{+\mu n}
-\bigl(\bigl(\me P^{\mu}_N\bigr)^0\bigr)^{-1}\bigl(\me T^{\mu
}_N\bigr)^0 \me a^{\mu}_N\bigr)_ke^{i\lambda
k}.
\]
This constrained optimization problem is equivalent to the
unconstrained optimization problem
\begin{equation}\label{zad_bezum_extr_i}
\varDelta_{\mathcal{D}}(f,g)=\widetilde{\varDelta}(f,g)+ \varDelta
(f,g|\mathcal
{D}_f\times
\mathcal{D}_g)\to\inf,
\end{equation}
where $\varDelta(f,g|\mathcal{D}_f\times
\mathcal{D}_g)$ is the indicator function of the set
$\mathcal{D}_f\times\mathcal{D}_g$: $\varDelta(f,g|\mathcal
{D}_f\times
\mathcal{D}_g)=0$ if $(f;g)\in\mathcal{D}_f\times
\mathcal{D}_g$ and $\varDelta(f,g|\mathcal{D}_f\times
\mathcal{D}_g)=+\infty$ if $(f;g)\notin\mathcal{D}_f\times
\mathcal{D}_g$.

\noindent A solution $(f^0,g^0)$ to the unconstrained optimization problem
is determined by the condition
$0\in\partial\varDelta_{\mathcal{D}}(f^0,g^0)$, which is a
necessary and
sufficient condition that the pair $(f^0,g^0)$ belongs to the set of
minimums of the convex functional $\varDelta_{\mathcal{D}}(f,g)$
\cite
{Tihomirov,Pshenychn,Rockafellar}. By $\partial\varDelta_{\mathcal
{D}}(f,g)$ we denote the subdifferential of the functional $\varDelta
_{\mathcal{D}}(f,g)$ at the point $(f,g)=(f^0,g^0)$, that is, the set
of all linear continuous functionals $\varLambda$ on the space
$L_1\times
L_1$ that satisfy the inequality
\[
\varDelta_{\mathcal{D}}(f,g)-\varDelta_{\mathcal{D}}\bigl
(f^0,g^0\bigr)\geq\varLambda
\bigl((f,g)-\bigl(f^0,g^0\bigr)\bigr),\quad(f,g)\in\mathcal{D}.
\]

In the case of estimating the cointegrated sequences, we have the
following optimization problem of finding the least favorable spectral
densities:
\begin{equation}\label{zad_bezum_extr_i_coi}
\varDelta_{\mathcal{D}}(f,p)=\widetilde{\varDelta}(f,p)+\varDelta
(f,p|\mathcal
{D}_f\times
\mathcal{D}_p)\to\inf,
\end{equation}
\begin{align*}
&\widetilde{\varDelta}(f,p)\\
&\quad =\varDelta\bigl(h^{\beta}_{\mu}\bigl(f^0,p^0\bigr);f,p\bigr)\\
&\quad =\frac{ 1}{2\pi}\ip\frac{|A_N(e^{i\lambda})(1-e^{i\lambda\mu})^{n}\beta^2f^0(\lambda)-\lambda^{2n}C_{\mu,N}^{\beta,0}(e^{i\lambda})|^2}{\lambda^{2n}|1-e^{i\lambda\mu}|^{2n}(p^0(\lambda))^2}p(\lambda)d\lambda\\
&\qquad -\frac{ \beta^2}{2\pi}\ip\frac{|A_N(e^{i\lambda})(1-e^{i\lambda\mu})^{n}\beta^2f^0(\lambda)-\lambda^{2n}C_{\mu,N}^{\beta,0}(e^{i\lambda})|^2}{\lambda^{2n}|1-e^{i\lambda\mu}|^{2n}l(p^0(\lambda))^2}f(\lambda)d\lambda\\
&\qquad \,{+}\,\frac{ 1}{2\pi}\!\ip\!\frac{|A_N(e^{i\lambda})(1-e^{i\lambda\mu})^{n}[p^0(\lambda)-\beta^2f^0(\lambda)]_+ +\lambda^{2n}C_{\mu,N}^{\beta,0}(e^{i\lambda})|^2}{\lambda^{2n}|1-e^{i\lambda\mu}|^{2n}(p^0(\lambda))^2}f(\lambda)d\lambda
\end{align*}
\[
C^{\beta,0}_{\mu,N}\bigl(e^{i\lambda}\bigr)=\sum_{k=0}^{N+\mu n}
\bigl(\bigl(\bigl(\me
P^{\mu,\beta}_N\bigr)^0\bigr)^{-1}\bigl(\bigl[D_N^{\mu}\me a_N\bigr]_{+\mu n}
-\bigl(\me T^{\mu,\beta}_N\bigr)^0 \me a^{\mu}_N\bigr)\bigr
)_ke^{i\lambda
k}.
\]
A solution $(f^0,p^0)$ to this optimization problem is characterized by
the condition
$0\in\partial\varDelta_{\mathcal{D}}(f^0,p^0)$.

The derived representations of the linear functionals $\varDelta
(h_{\mu
}(f^0,g^0);f,g)$ and\break $\varDelta(h^{\beta}_{\mu}(f^0,p^0);f,p)$
allow us
to calculate derivatives and subdifferentials in the space $L_1\times
L_1$. Therefore, the complexity of the optimization problems (\ref
{zad_bezum_extr_i}) and (\ref{zad_bezum_extr_i_coi}) is determined by
the complexity of calculation of the subdifferentials of the indicator
functions $\varDelta(f,g|\mathcal{D}_f\times
\mathcal{D}_g)$ and $\varDelta(f,p|\mathcal{D}_f\times
\mathcal{D}_p)$ of the sets $\mathcal{D}_f\times
\mathcal{D}_g$ and $\mathcal{D}_f\times
\mathcal{D}_p$.

\section{The least favorable spectral densities in the class $\md
D^-_{0,f}\times\md D^-_{0,g}$}

Consider the problem of minimax-robust estimation of the functional
$A_N\xi$ of unknown values of the sequence with stationary increments
$\xi(m)$ based on observations of the sequence $\xi(m)+\eta(m)$ at the
points $m\in\mr Z\setminus\{0,1,2,\ldots,N\}$ for the set of admissible
spectral densities
$\md D=\md D^-_{0,f}\times\md D^-_{0,g}$, where
\[
\md D^-_{0,f}=\ld\{f(\lambda)\bigg|\frac{1}{2\pi}\ip
\frac{1}{f(\lambda)}d\lambda\geq P_1\rd\},
\]
\[
\md
D^-_{0,g}=\ld\{g(\lambda)\bigg| \frac{1}{2\pi}\ip\frac
{1}{g(\lambda
)}d\lambda\geq
P_2\rd\}.
\]
If the spectral densities $f^0\in\md D^-_{0,f}$,
$g^0\in\md D^-_{0,g}$ and the functions
\begin{equation}
h_{\mu,f}\bigl(f^0,g^0\bigr)=\frac{|A_N(e^{i\lambda
})(1-e^{i\lambda
\mu})^{n}\lambda^{2n}g^0(\lambda)
+\lambda^{2n}C^{\mu,0}_{N}(e^{i\lambda})|}
{|\lambda
|^{n}|1-e^{i\lambda\mu}|^{n}p^0(\lambda)},\label{hf_i_st.n_d}
\end{equation}
\begin{equation}
h_{\mu,g}\bigl(f^0,g^0\bigr)=
\frac{|A_N(e^{i\lambda})(1-e^{i\lambda\mu})^{n}f^0(\lambda)- \lambda^{2n}C^{\mu,0}_{N}(e^{i\lambda})|}
{|1-e^{i\lambda\mu}|^{n}p^0(\lambda)}\label{hg_i_st.n_d},
\end{equation}
where $p^0(\lambda)=f^0(\lambda)+\lambda^{2n}g^0(\lambda)$, are
bounded, then the linear functional\break
$\varDelta(h_{\mu}(f^0,g^0);f,g)$ is continuous and bounded in the
space $
L_1\times L_1$. The condition
$0\in\partial\varDelta_{\md D}(f^0,g^0)$ implies that the spectral densities
$f^0\in\md D^-_{0,f}$ and $g^0\in\md D^-_{0,g}$ are determined by the
relations
\begin{align}
&|\lambda|^{n}f^0(\lambda)\big|A_N\bigl(e^{i\lambda}\bigr)\bigl(1-e^{i\lambda\mu}\bigr)^{n}g^0(\lambda)+C^{\mu,0}_{N}\bigl(e^{i\lambda}\bigr)\big|\notag\\
&\quad =\alpha_1\big|1-e^{i\lambda\mu}\big|^{n}\bigl(f^0(\lambda)+\lambda^{2n}g^0(\lambda)\bigr),\label{D1 rivn1_i_st.n_d}
\end{align}
\vspace*{-18pt}%
\begin{align}
&g^0(\lambda)\big|A_N\bigl(e^{i\lambda}\bigr)\bigl(1-e^{i\lambda
\mu}\bigr)^{n}f^0(\lambda)-\lambda^{2n}C^{\mu,0}_{N}\bigl
(e^{i\lambda}\bigr)\big|\notag
\\
&\quad =\alpha_2\big|1-e^{i\lambda\mu}\big|^{n}\bigl(f^0(\lambda
)+\lambda^{2n}g^0(\lambda
)\bigr),\label{D1
rivn2_i_st.n_d}
\end{align}
where the constants $\alpha_1\geq0$, $\alpha_2\geq0$ with
$\alpha_1\neq0$ if $\ip(f^0(\lambda))^{-1}d\lambda=2\pi P_1$ and
$\alpha_2\neq0$ if
$\ip(g^0(\lambda))^{-1}d\lambda=2\pi P_2$.

The derived statements allow us to formulate the following theorems.

\begin{thm} Suppose that the spectral densities $f^0(\lambda)\in\md
D^-_{0,f}$ and
$g^0(\lambda)\in\md D^-_{0,g}$ satisfy the minimality condition
\eqref
{umova11_i_st.n_d} and the functions $h_{\mu,f}(f^0,g^0)$ and
$h_{\mu,g}(f^0,g^0)$ calculated by formulas
\eqref{hf_i_st.n_d} and \eqref{hg_i_st.n_d}
are bounded. The spectral densities
$f^0(\lambda)$ and $g^0(\lambda)$ determined by Eqs.~\eqref{D1
rivn1_i_st.n_d}) and \eqref{D1 rivn2_i_st.n_d}
are the least favorable densities in the class $\md D=\md D^-_{0,f}
\times\md D^-_{0,g}$ for the linear
interpolation of the functional $A_N\xi$
if they give a solution to the constrained optimization problem~\eqref
{minimax1_i_st.n_d}. The function $h_{\mu}(f^0,g^0)$ calculated by
formula \eqref{spectr A_i_st.n_d}
is the minimax-robust spectral characteristic of the optimal estimate
of the functional
$A_N\xi$.
\end{thm}

\begin{thm}
Suppose that the spectral density $f(\lambda)$ \textup{(}or
$g(\lambda
))$ is known, the spectral density
$g^0(\lambda)\in\md D^-_{0,g}$ $(f^0(\lambda)\in\md D^-_{0,f})$, and
they satisfy the minimality condition~\eqref{umova11_i_st.n_d}. Suppose
also that the function $h_{\mu,g}(f ,g^0)$ $(h_{\mu,f}(f^0 ,g))$ is bounded.
Then the spectral density\vspace*{-3pt}
\[
g^0(\lambda)=f(\lambda)\ld[\frac{1}{\alpha_2|1-e^{i\lambda
\mu}|^{n}}\big
|A_N\bigl(e^{i\lambda}\bigr)\bigl(1-e^{i\lambda
\mu}\bigr)^{n}f(\lambda)-C^{\mu,0}_{N}\bigl(e^{i\lambda}\bigr
)\big|-\lambda^{2n}\rd]_+^{-1}
\]
or\vspace*{-3pt}
\[
f^0(\lambda)=\lambda^{2n}g(\lambda)\ld[\frac{|\lambda
|^{n}}{\alpha
_1|1-e^{i\lambda\mu}|^{n}}\big|A_N\bigl(e^{i\lambda}\bigr
)\bigl(1-e^{i\lambda
\mu}\bigr)^{n}g(\lambda)+C^{\mu,0}_{N}\bigl(e^{i\lambda}\bigr
)\big|-1\rd]_+^{-1}
\]\eject
\noindent is the least
favorable in the class
$\md D^-_{0,g}$ \textup{(}or $\md D^-_{0,f})$ for the linear
interpolation of the functional $A_N\xi$
if the functions $f(\lambda)+\lambda^{2n}g^0(\lambda)$, $g^0(\lambda)$
\textup{(}or $f^0(\lambda)+\lambda^{2n}g(\lambda))$ give a solution to
the constrained optimization problem
\eqref{minimax1_i_st.n_d}. The function $h_{\mu}(f ,g^0)$ \textup{(}or
$h_{\mu}(f^0,g))$ calculated by formula \eqref{spectr A_i_st.n_d} is
the minimax-robust spectral characteristic of the optimal estimate of
the functional
$A_N\xi$.
\end{thm}

Consider the problem of minimax-robust estimation of the functional
$A_N\xi$ of unknown values of the sequence $\xi(m)$, cointegrated with
the sequence $\zeta(m)$, based on observations of the sequence $\zeta
(m)$ at the points $m\in\mr Z\setminus\{0,1,2,\ldots,N\}$. Suppose that
the stochastic sequences $\xi(m)$ and $\zeta(m)-\beta\xi(m)$ are
uncorrelated.
The least favorable spectral densities in the class $\md D^0_f\times
\md
D^0_p$, where
\[
\md D^-_{0,f}\,{=}\,\ld\{f(\lambda)\bigg|\frac{1}{2\pi}\ip
\frac{1}{f(\lambda)}d\lambda\geq P_1\rd\},\quad\md
D^-_{0,p}\,{=}\,\ld\{p(\lambda)\bigg| \frac{1}{2\pi}\ip\frac
{1}{p(\lambda
)}d\lambda\geq
P_2\rd\},
\]
are determined by the condition $0\in\partial\varDelta_{\mathcal
{D}}(f^0,p^0)$, which implies
the following relations for determining the least favorable spectral
densities $f^0\in\md D^0_f$ and $p^0\in\md D^0_p$:
\begin{equation}
\big|A_N\bigl(e^{i\lambda}\bigr)\bigl(1-e^{i\lambda\mu}\bigr
)^{n}\beta^2f^0(\lambda)
-\lambda^{2n}C_{\mu,N}^{\beta,0}\bigl(e^{i\lambda}\bigr)\big
|=\alpha_2
|\lambda|^{n}\big|1-e^{i\lambda\mu}\big|^{n},\label{D1
rivn1_coi_i_st.n_d}
\end{equation}
\vspace*{-12pt}%
\begin{align}
&f^0(\lambda)\big|A_N\bigl(e^{i\lambda}\bigr)\bigl(1-e^{i\lambda
\mu}\bigr)^{n}\bigl[p^0(\lambda)-\beta^2f^0(\lambda)\bigr]_+
+\lambda^{2n}C_{\mu,N}^{\beta,0}\bigl(e^{i\lambda}\bigr)\big|
\notag\\
&\quad =|\lambda|^{n}\big|1-e^{i\lambda\mu}\big|^{n}\bigl(\alpha
_1p^0(\lambda)+\alpha
_2|\beta|f^0(\lambda)\bigr),\label{D1
rivn2_coi_i_st.n_d}
\end{align}
where the constants $\alpha_1\geq0$, $\alpha_2\geq0$ with $\alpha
_1\neq
0$ if $\ip(f^0(\lambda))^{-1}d\lambda=2\pi P_1$ and $\alpha_2\neq
0$ if
$\ip(p^0(\lambda))^{-1}d\lambda=2\pi P_2$.

\begin{thm} Suppose that the spectral density
$p^0(\lambda)\in\md D^-_{0,p}$ satisfies the minimality condition
\eqref
{umova111_i_st.n_d} and the functions $h_{\mu,f}(f^0,g^0)$ and
$h_{\mu,g}(f^0,g^0)$, calculated by formulas
\eqref{hf_i_st.n_d} and \eqref{hg_i_st.n_d}, are bounded for
$g(\lambda
):=\lambda^{-2n}(p(\lambda)-\beta^2f(\lambda))$. The spectral densities
$f^0(\lambda)$ and $p^0(\lambda)$ determined by Eqs.~\eqref{D1
rivn1_coi_i_st.n_d} and \eqref{D1 rivn2_coi_i_st.n_d} are the least
favorable in the class $\md D=\md D^-_{0,f} \times\md D^-_{0,p}$ for
the linear
interpolation of the functional $A_N\xi$ based on observations of the
stochastic sequence $\zeta(m)$, which is cointegrated with $\xi(m)$ and
such that the stochastic sequences $\xi(m)$ and $\zeta(m)-\beta\xi(m)$
are uncorrelated,
if these densities determine a solution to constrained optimization
problem \eqref{minimax1_i_st.n_d} for $g^0(\lambda):=\lambda
^{-2n}(p^0(\lambda)-\beta^2f^0(\lambda))$. The function $h_{\mu
}(f^0,p^0)$, calculated by formula \eqref{spectr A_N_coi_i_st.n_d},
is the minimax-robust spectral characteristic of the optimal estimate
of the functional
$A_N\xi$.
\end{thm}

\section{The least favorable spectral densities in the class $\md
D=\md
D_{2\varepsilon_1}\times\md D_{1\varepsilon_2}$}

Consider the problem of minimax-robust interpolation of the functional
$A_N\xi$ based on observations of the sequence $\xi(m)+\eta(m)$ at the
points of $m\in\mr Z\setminus\{0,1,2,\ldots,N\}$ in the case where the
spectral densities $f(\lambda)$ and $g(\lambda)$ belong to the set
$\md
D=\md D_{2\varepsilon_1}\times\md
D_{1\varepsilon_2}$, where
\[
\md D_{2\varepsilon_1}=\ld\{f(\lambda)\bigg|\frac{1}{2\pi}\ip
\big|f(\lambda
)-f_1(\lambda)\big|^2 d\lambda\leq\varepsilon_1\rd\},
\]
\[
\md
D_{1\varepsilon_2}=\ld\{g(\lambda)\bigg|\frac{1}{2\pi}\ip
\big|g(\lambda
)-g_1(\lambda)\big| d\lambda\leq\varepsilon_2\rd\}
\]
are $\varepsilon$-neighborhoods of the given spectral densities
$f_1(\lambda)$ and $g_1(\lambda)$ in the spaces $L_2$ and $L_1$, respectively.

Suppose that the spectral densities $f_1(\lambda)$ and $g_1(\lambda)$
are bounded and
the functions $h_{\mu,f}(f^0,g^0)$ and $h_{\mu,g}(f^0,g^0)$ calculated
by formulas (\ref{hf_i_st.n_d}) and
(\ref{hg_i_st.n_d}) with spectral densities $f^0\in
\md D_{2\varepsilon
_1}$ and $g^0\in\md D_{1\varepsilon_2}$ are bounded as well.
The condition
$0\in\partial\varDelta_{\md D}(f^0,g^0)$ implies the following relations
for determining the least favorable spectral densities:
\begin{align}
&\big|A_N\bigl(e^{i\lambda}\bigr)\bigl(1-e^{i\lambda
\mu}\bigr)^{n}\lambda^{2n}g^0(\lambda)+\lambda^{2n}C^{\mu
,0}_{N}\big(e^{i\lambda
}\big)\big|^2
\notag\\
&\quad =\alpha_1|\lambda|^{2n}\big|1-e^{i\lambda\mu}\big|^{2n}\bigl
(f^0(\lambda
)-f_1(\lambda)\bigr)\bigl(f^0(\lambda)+\lambda^{2n}g^0(\lambda
)\bigr)^2
,\label{D2 rivn1_i_st.n_d}
\end{align}
\vspace*{-18pt}%
\begin{align}
&\big|A_N\bigl(e^{i\lambda}\bigr)\bigl(1-e^{i\lambda
\mu}\bigr)^{n}f^0(\lambda)-\lambda^{2n}C^{\mu,0}_{N}\bigl
(e^{i\lambda}\bigr)\big|^2
\notag\\
&\quad =\alpha_2\gamma(\lambda)\big|1-e^{i\lambda\mu}\big|^{2n}\bigl
(f^0(\lambda)+\lambda
^{2n}g^0(\lambda)\bigr)^2,\label{D2 rivn2_i_st.n_d}
\end{align}
where the function
$|\gamma(\lambda)|\leq1$ and $\gamma(\lambda)=\text
{sign}(g(\lambda
)-g_1(\lambda)) $ if $g(\lambda)\neq g_1(\lambda)$; $\alpha_1$,~$\alpha
_2$ are two constants to be found using the equations
\begin{equation}\label{D2 rivn4_i_st.n_d}
\frac{1}{2\pi}\ip\big|f^0(\lambda)-f_1(\lambda)\big|^2 d\lambda
=\varepsilon
_1,\qquad\frac{1}{2\pi}\ip\big|g^0(\lambda)-g_1(\lambda)\big|
d\lambda
=\varepsilon_2.
\end{equation}

Now we can present the following theorems, which describe the least
favorable spectral densities in the class $\md D=\md D_{2\varepsilon
_1}\times\md
D_{1\varepsilon_2}$.

\begin{thm}
Suppose that the spectral densities $f^0(\lambda)\in\md
D_{2\varepsilon
_1}$ and $g^0(\lambda)\in\md
D_{1\varepsilon_2}$ satisfy the minimality condition \eqref
{umova11_i_st.n_d}, the functions $h_{\mu,f}(f^0,g^0)$ and
$h_{\mu,g}(f^0,g^0)$, calculated by formulas \eqref{hf_i_st.n_d} and
\eqref{hg_i_st.n_d}, are bounded. The spectral densities
$f^0(\lambda)$ and $g^0(\lambda)$ determined by equations
\eqref{D2 rivn1_i_st.n_d}--\eqref{D2 rivn4_i_st.n_d} are the least
favorable spectral densities in the class $\md D=\md D_{2\varepsilon_1}
\times\md D_{1\varepsilon_2}$ for the linear
interpolation of the functional $A_N\xi$
if they give a solution to constrained optimization problem \eqref
{minimax1_i_st.n_d}. The function $h_{\mu}(f^0,g^0)$, calculated by
formula \eqref{spectr A_i_st.n_d}
is the minimax-robust spectral characteristic of the optimal estimate
of the functional
$A_N\xi$.
\end{thm}

\begin{thm}
Suppose that the spectral density $f(\lambda)$ is known, the spectral density
$g^0(\lambda)\in\md D_{1\varepsilon_2}$, and they satisfy the
minimality condition \eqref{umova11_i_st.n_d}. Suppose also that the
function $h_{\mu,g}(f ,g^0)$ calculated by formula \eqref{hg_i_st.n_d}
is bounded.
Then the spectral density
\[
g^0(\lambda)=\max\big\{g_1(\lambda),\lambda^{-2n}f_2(\lambda)\big
\},
\]
\[
f_2(\lambda)=\alpha_2^{-1}\big|1-e^{i\lambda\mu}\big|^{-n}\ld
|A_N\bigl(e^{i\lambda
}\bigr)\big(1-e^{i\lambda
\mu}\big)^{n}f(\lambda)-\lambda^{2n}C^{\mu,0}_{N}\bigl(e^{i\lambda
}\bigr)\rd
|-f(\lambda),
\]
is the least
favorable in the class
$\md D_{1\varepsilon_2}$ for the linear
interpolation of the functional $A_N\xi$ if a pair $(f,g^0)$ provides a
solution to constrained optimization problem
\eqref{minimax1_i_st.n_d}. The function $h_{\mu}(f,g^0)$, calculated by
formula \eqref{spectr A_i_st.n_d} is the minimax-robust spectral
characteristic of the optimal estimate of the functional $A_N\xi$.
\end{thm}

\begin{thm}
Suppose that the spectral density $g(\lambda)$ is known, the spectral density
$f^0(\lambda)\in\md D_{2\varepsilon_1}$, and they satisfy the
minimality condition \eqref{umova11_i_st.n_d}. Suppose also that the
function $h_{\mu,f}(f^0 ,g)$, calculated by formula \eqref
{hf_i_st.n_d}, is bounded. The spectral density
$f^0(\lambda)$ determined by the equation
\begin{align*}
&\big|A_N\bigl(e^{i\lambda}\bigr)\bigl(1-e^{i\lambda
\mu}\bigr)^{n}\lambda^{2n}g(\lambda)+\lambda^{2n}C^{\mu
,0}_{N}\bigl(e^{i\lambda
}\bigr)\big|^2
\\
&\quad =\alpha_1|\lambda|^{2n}\big|1-e^{i\lambda\mu}\big|^{2n}\bigl
(f^0(\lambda
)-f_1(\lambda)\bigr)\bigl(f^0(\lambda)+\lambda^{2n}g(\lambda
)\bigr)^2
\end{align*}
and the condition $ \ip|f^0(\lambda)-f_1(\lambda)|^2 d\lambda=2\pi
\varepsilon_1$ is the least
favorable spectral density in the class $\md D_{2\varepsilon_1}$ for
the linear
interpolation of the functional $A_N\xi$ if a pair $(f^0,g)$ provides a solution to constrained optimization problem
\eqref{minimax1_i_st.n_d}. The function $h_{\mu}(f^0,g)$ calculated by
formula \eqref{spectr A_i_st.n_d}
is the minimax-robust spectral characteristic of the optimal estimate
of the functional
$A_N\xi$.
\end{thm}

Consider the problem of minimax-robust interpolation of the functional
$A_N\xi$ in the case of cointegrated sequences $\xi(m)$ and $\zeta(m)$
on the set of admissible spectral densities $\md D=\md D_{2\varepsilon
_1}\times\md
D_{1\varepsilon_2}$, where
\[
\md D_{2\varepsilon_1}=\ld\{f(\lambda)\bigg|\frac{1}{2\pi}\ip
\big|f(\lambda
)-f_1(\lambda)\big|^2 d\lambda\leq\varepsilon_1\rd\},
\]
\[
\md
D_{1\varepsilon_2}=\ld\{p(\lambda)\bigg|\frac{1}{2\pi}\ip
\big|p(\lambda
)-p_1(\lambda)\big| d\lambda\leq\varepsilon_2\rd\}.
\]
From the condition
$0\in\partial\varDelta_{\md D}(f^0,g^0)$ we obtain the following relations
that determine the least favorable spectral densities:
\begin{align}
&\big|A_N\bigl(e^{i\lambda}\bigr)\bigl(1-e^{i\lambda\mu}\bigr
)^{n}\beta^2f^0(\lambda)
-\lambda^{2n}C_{\mu,N}^{\beta,0}\bigl(e^{i\lambda}\bigr)\big|^2
\notag\\
&\quad =\alpha_2
\lambda^{2n}\gamma(\lambda)\big|1-e^{i\lambda\mu}\big|^{2n}\bigl
(p^0(\lambda
)\bigr)^2,\label{D2
rivn1_coi_i_st.n_d}
\end{align}
\vspace*{-12pt}%
\begin{align}
&\big|A_N\bigl(e^{i\lambda}\bigr)\bigl(1-e^{i\lambda
\mu}\bigr)^{n}\bigl[p^0(\lambda)-\beta^2f^0(\lambda)\bigr]_+
+\lambda^{2n}C_{\mu,N}^{\beta,0}\bigl(e^{i\lambda}\bigr)\big|^2
\notag\\
&\quad =\lambda^{2n}\big|1-e^{i\lambda\mu}\big|^{2n}\bigl(p^0(\lambda
)\bigr)^2\bigl(\alpha
_1\bigl(f^0(\lambda)-f_1(\lambda)\bigr)+\alpha_2\beta^2
\gamma(\lambda)\bigr),\label{D2
rivn2_coi_i_st.n_d}
\end{align}
where the function
$|\gamma(\lambda)|\leq1$ and $\gamma(\lambda)=\text
{sign}(p(\lambda
)-p_1(\lambda)) $ if $p(\lambda)\neq p_1(\lambda)$; $\alpha_1$,~$\alpha
_2$ are two constants that can be found from the equations
\begin{equation}\label{D2 rivn4_coi_i_st.n_d}
\frac{1}{2\pi}\ip\big|f^0(\lambda)-f_1(\lambda)\big|^2 d\lambda
=\varepsilon
_1,\qquad\frac{1}{2\pi}\ip\big|p^0(\lambda)-p_1(\lambda)\big|
d\lambda
=\varepsilon_2.
\end{equation}

Thus, we have the following theorem.

\begin{thm} Suppose that the spectral density
$p^0(\lambda)\in\md D_{1\varepsilon_2}$ satisfies the minimality
condition \eqref{umova111_i_st.n_d} and the functions $h_{\mu
,f}(f^0,g^0)$ and
$h_{\mu,g}(f^0,g^0)$, calculated by formulas \eqref{hf_i_st.n_d} and
\eqref{hg_i_st.n_d}, are bounded for $g(\lambda):=\lambda
^{-2n}(p(\lambda)-\beta^2f(\lambda))$. Then the least favorable
spectral densities for
the linear
interpolation of the functional $A_N\xi$ based on observations of the
stochastic sequence $\zeta(m)$, which is cointegrated with $\xi(m)$ and
such that the stochastic sequences $\xi(m)$ and $\zeta(m)-\beta\xi(m)$
are uncorrelated, are the spectral densities
$f^0(\lambda)$ and $p^0(\lambda)$ determined by
Eqs.~\eqref{D2 rivn1_coi_i_st.n_d}--\eqref{D2 rivn4_coi_i_st.n_d} and
provide a solution to constrained optimization problem \eqref
{minimax1_i_st.n_d} for $g^0(\lambda):=\lambda^{-2n}(p^0(\lambda
)-\beta
^2f^0(\lambda))$. The function $h_{\mu}(f^0,p^0)$ calculated by formula
\eqref{spectr A_N_coi_i_st.n_d}
is the minimax-robust spectral characteristic of the optimal estimate
of the functional
$A_N\xi$.
\end{thm}
\section{Conclusions}

In the article, the problem of the mean-square optimal linear
estimation of the functional $A_N\xi=\sum_{k=0}^{N}a(k)\xi(k)$, which
depends of unknown values of the sequence $\xi(m)$ with $n$th
stationary increments based on observations of the sequence $\xi
(m)+\eta
(m)$ at the points $m\in\mr Z\setminus\{0,1,2,\ldots,N\}$, is
considered in the case of observations with the stationary noise $\eta
(m)$ uncorrelated with $\xi(m)$. The classical and minimax-robust
methods of interpolation are applied in the case of spectral certainty
and in the case spectral uncertainty.
Particularly, in the case of spectral certainty, formulas for
calculating the spectral characteristics and the value of the
mean-square error of the optimal estimate are found. The derived
results are applied to interpolation problem for a class of
cointegrated sequences. In the case spectral uncertainty, where
spectral densities are not known exactly, whereas some
sets of admissible spectral densities are given, formulas that
determine the least
favorable spectral densities and the minimax-robust spectral
characteristics are derived for some special sets of admissible
spectral densities.
\section*{Acknowledgments}
The authors would like to thank the referee for careful
reading of the article and giving constructive suggestions.

%
\end{document}